# THE MINIMUM 3-COVERING ENERGY OF COMPLETE GRAPHS

PAUL AUGUST WINTER

DEPARTMENT OF MATHEMATICS, UNIVERSITY OF KWAZULU NATAL, HOWARD COLLEDGE, GLENWOOD, DURBAN, SOUTH AFRICA, 4041

email: winter@ukzn.ac.za

ABSTRACT: In this paper we introduce a different kind of graph energy, the minimum 3-covering energy of a graph, and determine the minimum 3-covering energy of complete graphs.

## 1. INTRODUCTION

The Huckel Molecular Orbital theory provided the motivation for the idea of the *energy* of a graph – the sum of the absolute values of the eigenvalues associated with the graph (see [1]). This resulted in the idea of the *minimum 2-covering energy* of a graph in [1]. This idea is generalized to the minimum 3-covering energy of a graph in this paper.

All graphs which we shall consider will be finite, simple, loopless and undirected. Let G be such a graph of order n with vertex set $\{v_1, v_2, ..., v_n\}$. A covering (2-covering) of a graph G is a set S of vertices of G of such that *every edge* of G has at least one vertex in S (see [1]) Since an edge is a path length 1 on 2 vertices ( a 2-path) we generalize this to a *3-covering* of a graph G as being set Q of vertices of G such that every path of G of length 2 (or 3-path) has *at least one vertex* in Q (see [2]). Any 3-covering set of G of minimum cardinality is called a *minimum 3-covering* of G.

## 2. THE MINIMUM 3-COVERING ENERGY OF A GRAPH

A minimum 3-covering matrix of G with a minimum 3-covering set Q of vertices is a matrix:

$$A_Q^3(G) = (a_{i,j})$$

where

$$a_{ij} = \begin{cases} 1 & \text{if } v_i v_j \in E(G) \\ 1 & \text{if } i = j \text{ and } v_i \in Q \\ 0 & \text{otherwise} \end{cases} (*)$$

The middle condition (*) is equivalent to loops of weight 1 being attached to the vertices of Q.

The characteristic polynomial of $A_Q^3(G)$ is then denoted by

$$f_n(G, \lambda) := \det(\lambda I - A_Q^3(G))$$

The *minimum 3-covering energy* is then defined as:

$$E_Q(G) = \sum_1^n |\lambda_i|$$

Where $\lambda_i$ (the *minimum 3-covering eigenvalues*) are the n real roots of the characteristic polynomial.

## 3. MOLECULAR STRUCTURES AND ENERGY

The minimum 2-covering energy of molecular structures involves the smallest set of atoms, such that every atom of the structure, is either in the set, or is connected (via bonding) *directly* to at least one vertex of the set. This is generalized to a minimum 3-covering energy of molecular structures, where the smallest set of atoms is considered, such that every path of 3 atoms has at least one atom in this set, so that no atom is more than a distance 2 from this set.

## 4. EXAMPLES

For example, consider the path $P_3$ on vertices $v_1, v_2, v_3$ where the first and the last listed vertices are end (pendant) vertices.

A minimum 3-covering is $\{v_1\}$ so that:

$$A_Q^3(P_3) = \begin{bmatrix} 1 & 1 & 0 \\ 1 & 0 & 1 \\ 0 & 1 & 0 \end{bmatrix} \quad Q = \{v_1\}$$

The characteristic polynomial is therefore:

$$\det(\lambda I - A_Q^3(P_3)) = \begin{vmatrix} \lambda-1 & -1 & 0 \\ -1 & \lambda & -1 \\ 0 & -1 & \lambda \end{vmatrix}$$

$$= \begin{vmatrix} \lambda-1 & 0 \\ -1 & -1 \end{vmatrix} + \lambda \begin{vmatrix} \lambda-1 & -1 \\ -1 & \lambda \end{vmatrix} = 1 - \lambda + \lambda(\lambda^2 - \lambda - 1) = \lambda^3 - \lambda^2 - 2\lambda + 1$$

The eigenvalues are (to 5 decimal places – online bluebit matrix calculator):

0.44504; -1,24698; 1,80134 so that the minimum 3-covering energy of this path is:

3,49396

If we take $\{v_2\}$ as another minimum 3-covering of the same path then:

$$A_Q^3(P_3) = \begin{bmatrix} 0 & 1 & 0 \\ 1 & 1 & 1 \\ 0 & 1 & 0 \end{bmatrix} \quad Q = \{v_2\}$$

The characteristic polynomial is therefore:

$$\det(\lambda I - A_Q^3(P_3)) = \begin{vmatrix} \lambda & -1 & 0 \\ -1 & \lambda-1 & -1 \\ 0 & -1 & \lambda \end{vmatrix}$$

$$= \begin{vmatrix} \lambda & 0 \\ -1 & -1 \end{vmatrix} + \lambda \begin{vmatrix} \lambda & -1 \\ -1 & \lambda-1 \end{vmatrix} = -\lambda + \lambda(\lambda^2 - \lambda - 1) = \lambda^3 - \lambda^2 - 2\lambda = \lambda(\lambda^2 - \lambda - 2)$$

$\lambda(\lambda - 2)(\lambda + 1)$.

The minimum 3-covering eigenvalues are therefore: 0,2 and -1 so that the minimum 3-covering energy of this path is 3.

The above examples illustrate that, in case of a path on 3 vertices, the minimum 3-covering energy of a graph depends on the choice of the 3-covering set.

However, consider the completer graph on 3 vertices, vertices labeled $v_1, v_2, v_3$. We take our 3-covering set as $v_1$.

$$A_Q^3(K_3) = \begin{bmatrix} 1 & 1 & 1 \\ 1 & 0 & 1 \\ 1 & 1 & 0 \end{bmatrix} \quad Q = \{v_1\}$$

The characteristic polynomial is therefore:

$$\det(\lambda I - A_Q^3(K_3)) = \begin{vmatrix} \lambda - 1 & -1 & -1 \\ -1 & \lambda & -1 \\ -1 & -1 & \lambda \end{vmatrix}$$

$$= \begin{vmatrix} \lambda - 1 & -1 & -1 \\ -1 & \lambda & -1 \\ 0 & -1-\lambda & \lambda+1 \end{vmatrix}$$

$$= (1+\lambda)\begin{vmatrix} \lambda - 1 & -1 \\ -1 & -1 \end{vmatrix} + (1+\lambda)\begin{vmatrix} \lambda - 1 & -1 \\ -1 & \lambda \end{vmatrix}$$

$$= (1+\lambda)(-\lambda) + (1+\lambda)(\lambda^2 - \lambda - 1) = (1+\lambda)(\lambda^2 - 2\lambda - 1)$$

So that eigenvalues are $-1, \dfrac{2 \pm \sqrt{8}}{2}$ so that the minimum 3-covering energy of the complete graph is $1 + \sqrt{8}$.

If we take our minimum 3-covering as $v_2$ instead of $v_1$ then we have:

$$A_Q^3(K_3) = \begin{bmatrix} 0 & 1 & 1 \\ 1 & 1 & 1 \\ 1 & 1 & 0 \end{bmatrix} \quad Q = \{v_2\}$$

the characteristic polynomial is:

$$\det(\lambda I - A_Q^3(P_4)) = \begin{vmatrix} \lambda & -1 & -1 \\ -1 & \lambda-1 & -1 \\ -1 & -1 & \lambda \end{vmatrix} = \begin{vmatrix} \lambda & -1 & -1 \\ -1 & \lambda-1 & -1 \\ 0 & -\lambda & \lambda+1 \end{vmatrix}$$

$$= \lambda \begin{vmatrix} \lambda & -1 \\ -1 & -1 \end{vmatrix} + (1+\lambda) \begin{vmatrix} \lambda & -1 \\ -1 & \lambda-1 \end{vmatrix}$$

$$= \lambda(-1-\lambda) + (1+\lambda)(\lambda^2 - \lambda - 1)$$

$$(1+\lambda)(\lambda^2 - 2\lambda - 1)$$

Which is the same as when Q was a different vertex.

## 5. THE MIMIMUM 3-COVERING ENERGY OF COMPLETE GRAPHS

Generally, the minimum 2-covering of a complete graph G on n vertices is any set of (n-1) vertices of G (see 1). The minimum 3-covering of a complete graph G on n vertices is any set of (n-2) vertices. Thus:

$$A_Q^3(K_n) = \begin{bmatrix} 1 & 1 & 1 & : & 1 & 1 \\ 1 & 1 & 1 & : & 1 & 1 \\ 1 & 1 & : & 1 & 1 & : \\ : & 1 & 1 & 1 & : & 1 \\ 1 & : & 1 & 1 & 0 & 1 \\ 1 & 1 & : & 1 & 1 & 0 \end{bmatrix}_{n \times n}$$ ; Hence the characteristic equation:

$$\det(\lambda I - A_Q^3(K_n)) = \det\begin{bmatrix} \lambda-1 & -1 & -1 & : & -1 & -1 \\ -1 & \lambda-1 & -1 & : & -1 & -1 \\ -1 & -1 & : & -1 & -1 & : \\ : & -1 & -1 & \lambda-1 & : & -1 \\ -1 & : & -1 & -1 & \lambda & -1 \\ -1 & -1 & : & -1 & -1 & \lambda \end{bmatrix}$$

Subtracting the last row from the second to last row:

$$= \det\begin{bmatrix} \lambda-1 & -1 & -1 & : & -1 & -1 \\ -1 & \lambda-1 & -1 & : & -1 & -1 \\ -1 & -1 & : & -1 & -1 & : \\ : & -1 & -1 & \lambda-1 & : & -1 \\ -1 & : & -1 & -1 & \lambda & -1 \\ 0 & 0 & : & 0 & -1-\lambda & \lambda+1 \end{bmatrix}_{nxn}$$

Expanding this determinant using the last row yields:

$$(1+\lambda)\begin{vmatrix} \lambda-1 & -1 & : & -1 & -1 \\ -1 & \lambda-1 & -1 & : & -1 \\ -1 & -1 & : & -1 & : \\ : & : & -1 & \lambda-1 & -1 \\ -1 & -1 & -1 & : & -1 \end{vmatrix}_{((n-1)n-1))}$$

$$(\lambda+1)\begin{vmatrix} \lambda-1 & -1 & : & -1 & -1 \\ -1 & \lambda-1 & -1 & : & -1 \\ -1 & -1 & : & -1 & : \\ -1 & : & -1 & \lambda-1 & -1 \\ -1 & -1 & : & -1 & \lambda \end{vmatrix}_{(n-1)(n-1)}$$

The first determinant is $-\lambda^{n-2}$, and we subtract the last row of the second determinant from the second to last row:

$$= -(1+\lambda)\lambda^{n-2} + (\lambda+1)\begin{vmatrix} \lambda-1 & -1 & : & -1 & -1 \\ -1 & \lambda-1 & -1 & : & -1 \\ -1 & -1 & : & -1 & : \\ -1 & : & -1 & \lambda-1 & -1 \\ 0 & 0 & : & -\lambda & \lambda+1 \end{vmatrix}_{(n-1)(n-1)}$$

Expandint the determinat using the last row yields:

$$= -(1+\lambda)\lambda^{n-2} + (\lambda+1)\lambda\begin{vmatrix} \lambda-1 & -1 & : & -1 & -1 \\ -1 & : & -1 & : & -1 \\ -1 & -1 & \lambda-1 & -1 & : \\ -1 & : & -1 & \lambda-1 & -1 \\ -1 & -1 & -1 & -1 & -1 \end{vmatrix}_{(n-2)(n-2)}$$

$$+ (\lambda+1)(1+\lambda)\begin{vmatrix} \lambda-1 & -1 & : & -1 & -1 \\ -1 & \lambda-1 & -1 & : & -1 \\ -1 & -1 & : & -1 & : \\ -1 & : & -1 & \lambda-1 & -1 \\ -1 & -1 & : & -1 & \lambda-1 \end{vmatrix}_{(n-2)(n-2)}$$

The first determinant is $-\lambda^{n-3}$ while the second determinant yields $\lambda^{n-3}(\lambda-(n-2))$ so that we have:

$$-(1+\lambda)\lambda^{n-2} - (1+\lambda)\lambda^{n-2} + (1+\lambda)^2(\lambda^{n-3}(\lambda-(n-2)))$$

$$-2(1+\lambda)\lambda^{n-2} + (1+\lambda)^2\lambda^{n-2} - \lambda^{n-3}(1+\lambda)^2(n-2)$$
$$= (1+\lambda)\lambda^{n-3}(-2\lambda + (1+\lambda)\lambda - (1+\lambda)(n-2))$$
$$= (1+\lambda)\lambda^{n-3}(\lambda^2 - (n-1)\lambda - (n-2))$$

Eigenvalues are 0 (n-3 times), -1 and the conjugate pairs:

$$\frac{(n-1) \pm \sqrt{(n-1)^2 + (4n-8)}}{2} = \frac{(n-1) \pm \sqrt{n^2 + 2n - 7}}{2}$$

Thus we have the following theorem:

**THEOREM**

The minimum 3-covering energy of a complete graph on $n \geq 3$ vertices is:

$$1 + \sqrt{n^2 + 2n - 7}$$